\documentclass[12pt]{article}
\usepackage {amsmath, amssymb}
\def\ep{\varepsilon}

\begin{document}
\centerline{\bf\large LOG-LEVEL COMPARISON PRINCIPLE}
\centerline{\bf\large FOR SMALL BALL PROBABILITIES}
\bigskip

\centerline {\rm Alexander NAZAROV\footnote{Partially supported by
grant NSh.227.2008.1 and by RFBR grant No.07-01-00159
}}

\bigskip

\centerline {\it   St.Petersburg State University}

\bigskip

{ \it \footnotesize We prove a new variant of comparison principle for
logarithmic $L_2$-small ball probabilities of Gaussian processes. As an
application, we obtain logarithmic small ball asymptotics for some well-known
processes with smooth covariances.}

\bigskip


\bigskip

\section{Introduction}

The theory of small deviations of random functions is currently in intensive
development. In this paper we consider the most explored case of Gaussian
processes in $L_2$-norm.

Suppose we have a Gaussian random function $X(x)$, $x\in\Omega$, with zero
mean and covariance function $G_X(x,y)=EX(x)X(y)$ for $x,y\in\Omega$.
Let $\mu$ be a measure on $\Omega$. Set
\[
||X||_{\mu}=\Big(\int\limits_{\Omega} X^2(x)\,\mu(dx)\Big)^{1/2}.
\]

\noindent (if $\mu$ is the Lebesgue measure, the index $\mu$ will
be omitted). The problem is to define the behavior of ${\bf
P}\{||X||_{\mu}\leq \ep \}$ as $\ep \rightarrow 0$.

The study of small deviation problem was initiated by Sytaya \cite{S}
and continued by many authors. The history of the problem in 20th century is
described in reviews by Lifshits \cite{Lf} and by Li and Shao \cite{LS}.
Latest results can be found in \cite{site}.

According to the well-known Karhunen-Lo\`eve expansion, we
have in distribution
$$||X||_{\mu}^2 =\int\limits_{\Omega} X^2(x)
\mu(dx)\overset {d}{=} \sum_{n=1}^\infty \lambda_n\xi_n^2,\eqno (1.1)$$

\noindent where $\xi_n$, $n\in\mathbb N$, are independent standard
normal r.v.'s, and $\lambda_n>0$, $n\in\mathbb N$,
$\sum\limits_n\lambda_n <\infty$, are the eigenvalues of the
integral equation
$$ \lambda f(x)=\int\limits_{\Omega}
G_X(x,y)f(y)\mu(dy), \quad x\in \Omega.
\eqno(1.2) $$
Thus we need to study the asymptotic behavior as $\ep \to 0 $ of
${\bf P}\left\{\sum_{n=1}^\infty \lambda _n \xi_n^2\leq\ep ^2 \right\}$.
The answer heavily depends on the available information about the
eigenva\-lues $\lambda_n$. We underline that the explicit formulas for these
eigenvalues are known only for a limited number of Gaussian processes.

A quite important contribution to $L_2$-small ball problem was made by Li
\cite{Li} who established the so-called comparison theorem for sums (1.1). In
a slightly sharpened form, see \cite{GHLT}, it reads as follows.\medskip

{\bf Proposition 1}. { \it Let $(\lambda_n)$ and
$(\widetilde\lambda_n)$, $n\in\mathbb N$, be two positive summable
sequences. If the infinite product ${\cal P}=\prod_{n=1}^{\infty}
({\lambda_n}\big/{\widetilde\lambda_n})$ converges, then, as
$r\to0$,}
$${\bf P}\left\{\sum_{n=1}^\infty\widetilde\lambda_n\xi_n^2\le r\right\}
\sim {\cal P}^{\frac 12}\cdot {\bf
P}\left\{\sum_{n=1}^\infty\lambda_n\xi_n^2\le r\right\}.\eqno(1.3)
$$
Thus, if we know sufficiently sharp asymptotics of $\lambda_n$, we can, in
principle, calculate the asymptotics of ${\bf P}\{||X||_{\mu}\leq \ep \}$
up to a constant. Such asymptotics of eigenvalues can be obtained if the
function $X$ is a one-parameter Gaussian process (i.e. $\Omega$ is an interval)
and the covariance $G_X$ is the Green function of a boundary value problem
for ordinary differential operator. This approach was developed in \cite{NN1}
for the case of "separated"\ boundary conditions and in a recent work
\cite{Na1} in the general case. Note that if, in addition, the eigenfunctions
of the covariance kernel can be expressed in terms of elementary or special
functions, the sharp constants in the small ball asymptotics can be 
calculated explicitely by the complex variable methods \cite{Na2} (see also 
\cite{GHLT1}, \cite{Na1} and \cite{NP}).

In general case, we cannot expect to derive the exact asymptotics. So, we
need so-called {\bf logarithmic} asymptotics, that is the asymptotics of
$\ln {\mathbb P}\{||X||_{\mu}\leq \ep \}$ as $\ep\to0$. It was shown in
\cite{NN2} (see also \cite{KNN} and \cite{Na3}) that in some cases this
asymptotics is completely determined by the one-term asymptotics of
$\lambda_n$. This gives a logarithmic version of the comparison theorem.

However, this statement is not always satisfactory. Indeed, if $\lambda_n$
decrease too fast (for example, exponentially) then it is difficult to
calculate even the one-term asymptotics of eigenvalues in terms of covariance.
In this paper we present a new variant of the log-level comparison principle
based on the asymptotics of the {\bf counting function}
$${\cal N}(\lambda)=\#\{n\ :\ \lambda_n<\lambda\}.$$

{\bf Theorem 1}. { \it Let $(\lambda_n)$ and
$(\widetilde\lambda_n)$, $n\in\mathbb N$, be two positive summable
sequences with counting functions ${\cal N}(\lambda)$ and
$\widetilde{\cal N}(\lambda)$, respectively. Suppose ${\cal N}$
satisfies
$$\liminf\limits_{x\to0}\frac {\int\limits_0^{hx} {\cal N}(\lambda)\ d\lambda}
{\int\limits_0^{x}{\cal N}(\lambda)\ d\lambda}>1 \qquad
\mbox{for
any} \qquad h>1. \eqno(1.4)
$$
If
$${\cal N}(\lambda)\sim\widetilde{\cal N}(\lambda), \qquad \lambda\to 0,
\eqno(1.5)$$
then, as $r\to0$,}
$$\ln{\bf P}\left\{\sum_{n=1}^\infty\widetilde\lambda_n\xi_n^2
\le r \right\} \sim \ln{\bf P}\left\{ \sum_{n=1}^\infty
\lambda_n\xi_n^2\le r \right\}. \eqno(1.6)
$$

{\bf Remark 1}. For the power-type decreasing of $\lambda_n$ the
one-term asymptotics of ${\cal N}(\lambda)$ as $\lambda\to0$
provides the one-term asymptotics of $\lambda_n$. So, Proposition
2.1 \cite{NN2}, Theorem 4.2 \cite{KNN} and Theorem 4.2 \cite{Na3}
are particular cases of this statement. But for the super-power
decreasing of $\lambda_n$ the condition (1.5) is weaker than
$\lambda_n\sim\widetilde\lambda_n$.\medskip

{\bf Remark 2}. The assumption (1.4) is satisfied, for example,
for ${\cal N}(\lambda)$ regularly varying of order $p\in(-1,0)$ at zero,
see \cite{KNN}. However, if ${\cal N}(\lambda)=\frac
1{\lambda\ln^\sigma(\lambda)}$, $\sigma>1$, then the relation
(1.6) fails as it is shown in \cite[Proposition 4.4]{KNN}. Thus,
the assumption (1.4) cannot be removed.\medskip

The structure of the paper is as follows. In \S2 we prove Theorem 1. In \S3
we apply it to derive the logarithmic $L_2$-small ball asymptotics for the
Gaussian processes with smooth covariances. Another example of the process
with super-power decreasing of eigenvalues is considered in \cite{NSh}.

\section{Proof of Theorem 1}

We will use the following result which is a particular case of
Theorem 2 from \cite{Lf2}.\medskip

{\bf Proposition 2}. { \it Let $(\lambda_n)$, $n\in\mathbb N$, be
a positive summable sequence. Define, for $\ t,u\geq 0$, $f(t)=(1+2t)^{-1/2}$,
$$L(u)=\sum\limits_{n=1}^\infty \ln f(u\lambda_n). $$

\noindent Then, as $r \to 0$,
$$ P\left\{\sum_{n=1}^\infty
\lambda_n\xi_n^2\leq r\right\}\sim \frac {\exp (L(u)+ ur)} {\sqrt
{2\pi u^2 L''(u)}},\eqno (2.1) $$

\noindent where $u=u(r)$ is any function satisfying}
$$\lim_{r\to0}\frac {L'(u)+r}{\sqrt {L''(u)}}=0.\eqno (2.2) $$

We begin by the asymptotic analysis of $L'(u)$ as $u\to \infty$.
Clearly
$$L'(u)=-\sum\limits_{n=1}^\infty\frac {\lambda_n}
{1+2u\lambda_n}=\int\limits_0^A\frac {\lambda\ d{\cal N}(\lambda)}
{1+2u\lambda}$$ 
(here $A\ge\lambda_1$; without loss of
generality, one can suppose $A=1$). Note that this integral
converges since $(\lambda_n)$ is summable.

Integrating by parts twice we obtain
$$L'(u)=-\int\limits_0^1\frac {{\cal N}(\lambda)\
d\lambda}{(1+2u\lambda)^2}=- \frac {{\cal M}(1)}{(1+2u)^2}-
4u\int\limits_0^1\frac {{\cal M}(\lambda)\ d\lambda}
{(1+2u\lambda)^3}$$
where ${\cal
M}(\lambda)=\int\limits_0^\lambda{\cal N}(t)\ dt$.

Note that
$$u|L'(u)|\ge u\int\limits_0^{1/u}\frac {{\cal
N}(\lambda)\ d\lambda}{(1+2u\lambda)^2}=\int\limits_0^1\frac
{{\cal N}(\frac tu)\ dt}{(1+2t)^2}\ge\frac 13 {\cal N}(\frac 1u)
\to\infty, \qquad u\to\infty,$$

\noindent while, given $a<1$,
$$\frac {{\cal M}(1)}{(1+2u)^2}+4u\int\limits_a^1\frac {{\cal
M}(\lambda)\ d\lambda}{(1+2u\lambda)^3}=O(u^{-2}), \qquad
u\to\infty,$$

\noindent and hence
$$L'(u)=-4u\int\limits_0^a\frac {{\cal M}(\lambda)\ d\lambda}
{(1+2u\lambda)^3}\cdot (1+o(u^{-1}))=-4 \int\limits_0^{au}\frac
{{\cal M}(\frac tu)\ dt} {(1+2t)^3}\cdot(1+o(u^{-1})), \qquad
u\to\infty.
\eqno(2.3)$$

\noindent Using the Cauchy theorem and (1.4), we conclude that for any
given $h>1$ and sufficiently large $u$
$$\frac {L'(u)}{L'(hu)}=\frac {\int_0^u\frac {{\cal M}
(\frac tu)\ dt} {(1+2t)^3}\cdot(1+o(u^{-1}))} {\int_0^u\frac
{{\cal M}(\frac t{hu})\ dt} {(1+2t)^3}\cdot(1+o(u^{-1}))}=\frac
{{\cal M}(\frac {\widehat t} u)} {{\cal M}(\frac {\widehat t}
{hu})}\cdot(1+o(u^{-1}))\ge 1+\delta(h-1).
\eqno(2.4)$$

Let us define the function $\widetilde L(u)$ using the sequence
$(\widetilde\lambda_n)$ instead of $(\lambda_n)$. Due to (1.5),
given $\varepsilon>0$, we can find $a>0$ such that $\left|\frac
{\widetilde{\cal N}(t)}{{\cal N}(t)}-1\right| \le\varepsilon$ for
$t<a$, and (2.3) gives for sufficiently large $u$
$$\left|\frac {\widetilde L'(u)}{L'(u)}-1\right|=\left|\frac
{\int_0^a\frac {\widetilde{\cal M}(\lambda)\ d\lambda}
{(1+2u\lambda)^3}\cdot (1+o(u^{-1}))} {\int_0^a\frac {{\cal M}
(\lambda)\ d\lambda} {(1+2u\lambda)^3}\cdot (1+o(u^{-1}))}-1
\right|\le 2\varepsilon.
\eqno(2.5)$$

Now we observe that $L''(u)=4\int\limits_0^1\frac {\lambda{\cal
N}(\lambda)\ d\lambda}{(1+2u\lambda)^3}>0$, and therefore, the
equation $L'(u)+r=0$ has for sufficiently small $r$ the unique
solution $u(r)$ such that $u(r)\to\infty$ as $r\to0$. In a similar
way we define $\widetilde u(r)$ as the solution of the equation
$\widetilde L'(u)+r=0$.

Since $\widetilde L'(\widetilde u(r))\equiv L'(u(r))$, relations
(2.4) and (2.5) imply for arbitrary $\ep>0$ and for sufficiently small $r$
$$\delta\left(\frac{u(r)}{\widetilde u(r)}-1 \right)\le \frac
{L'(u(r))}{L'(\widetilde u(r))} -1=\frac {\widetilde L'(\widetilde
u(r))} {L'(\widetilde u(r))}-1\le 2\varepsilon,$$

\noindent and, similarly, $\delta\left(\frac{\widetilde u(r)}
{u(r)}-1 \right)\le 2\varepsilon$. Thus, $\widetilde u(r)\sim
u(r)$ as $r\to0$.

Since $u(r)$ trivially satisfies (2.2) we can apply formula (2.2)
with $u=u(r)$. Naturally, the replacement of $u(r)$ by its
one-term asymptotics $\widetilde u(r)$ breaks the relation (2.2).
Therefore, this replacement does not work to extract the explicit form 
of exact small ball asymptotics from (2.1). Fortunately, it works for the
logarithmic asymptotics. Namely, we have
$$L(u)=\frac 12 \int\limits_0^1\ln (1+2u\lambda)\ d{\cal N}
(\lambda)= -u\int\limits_0^1\frac {{\cal N}(\lambda)\
d\lambda}{1+2u\lambda},
$$

\noindent and formula (2.1) gives, as $r\to0$ and $u=u(r)$,
$$\ln P\left\{\sum_{n=1}^\infty \lambda_n\xi_n^2 \leq r \right\}
\sim L(u) +ur=L(u)-u L'(u)=- 2u^2\int\limits_0^1\frac {{\cal
N}(\lambda)\lambda\ d\lambda}{(1+2u\lambda)^2}.
$$

Note that
$$|L(u)|\ge 2u^2\int\limits_0^{1/u}\frac {{\cal N}(\lambda)\lambda\
d\lambda}{(1+2u\lambda)^2}=2\int\limits_0^1\frac {{\cal N}(\frac
tu)t\ dt}{(1+2t)^2} \to\infty, \qquad u\to\infty,$$

\noindent while, given $a<1$,
$$u^2\int\limits_a^1\frac {{\cal N}(\lambda)\lambda\ d\lambda}
{(1+2u\lambda)^2}=O(1),\qquad u\to\infty,$$

\noindent and hence, using $\widetilde u(r)\sim u(r)$, we obtain
$$\frac {\ln P\left\{\sum_{n=1}^\infty\widetilde \lambda_n \xi_n^2
\leq r \right\}} {\ln P\left\{\sum_{n=1}^\infty\lambda_n \xi_n^2
\leq r \right\}}\sim \frac {\widetilde u^2(r)\int_0^1 \frac
{\widetilde {\cal N}(\lambda)\lambda\ d\lambda}{(1+2\widetilde
u(r) \lambda)^2}} {u^2(r)\int_0^1 \frac {{\cal N}(\lambda)
\lambda\ d\lambda}{(1+2u(r)\lambda)^2}}\sim \frac {\int_0^a \frac
{\widetilde {\cal N}(\lambda)\lambda\ d\lambda}{(1+2u \lambda)^2}}
{\int_0^a \frac {{\cal N}(\lambda) \lambda\ d\lambda}
{(1+2u\lambda)^2}}.
$$

Due to (1.5), given $\varepsilon>0$, we can find $a>0$ such that
$\left|\frac {\widetilde{\cal N}(t)}{{\cal N}(t)}-1\right|
\le\varepsilon$ for $t<a$, and (1.6) follows.\hfill$\square$

\section{Examples}

The typical example of super-power eigenvalues decreasing is the
case where ${\cal N}(\lambda)$ is slowly varying at zero, see
\cite{Se}. We recall that this means
$$\lim\limits_{t\to0}\frac {{\cal N}(ct)}{{\cal N}(t)}=1
\qquad\mbox{for any}\quad c>0.\eqno(3.1)$$

{\bf Theorem 2}. { \it Let $(\lambda_n)$, $n\in\mathbb N$, be a
positive sequence with counting function ${\cal N}(\lambda)$.
Suppose that
$${\cal N}(\lambda)\sim \varphi(\lambda), \qquad \lambda\to 0,
\eqno(3.2)$$
where $\varphi$ is a function slowly varying at zero.
Then, as $r\to0$,
$$\ln{\bf P}\left\{\sum_{n=1}^\infty\lambda_n\xi_n^2\le r
\right\} \sim-\frac {\psi(\frac 1u)}2\equiv -\frac 12 \int\limits
_{1/u}^1 \varphi(z)\frac {dz}z, \eqno(3.3)
$$

\noindent where $u=u(r)$ satisfies the relation}
$$\frac {\varphi(\frac 1u)}{2u}\sim r, \qquad r\to 0.\eqno(3.4)$$

{\bf Remark 3}. According to \cite[Theorem 2.2]{KNN}, under
assumptions of Theorem 2 $\psi$ is also a slowly varying function.
Moreover, since $\varphi(z)\to\infty$ as $z\to0$, we have for
arbitrary large $A$ and $x>1/A$
$$\psi(x)\ge \int \limits_{x}^{Ax} \varphi(z)\frac{dz}z= \int
\limits_1^A \varphi(xt) \frac{dt}t\sim \varphi(x)\int\limits_1^A
\frac {dt}t= \ln(A)\cdot \varphi(x),
$$
which implies
$$\varphi(x)=o(\psi(x)),\qquad x\to 0,
\eqno(3.5)$$

{\bf Proof}. The relation (3.1) clearly implies
$$\lim\limits_{x\to0}\frac {\int\limits_0^{cx} {\cal N}(t)\ dt}
{\int\limits_0^{x}{\cal N}(t)\ dt}=c \qquad \mbox{for any} \qquad
c>0,
$$
and (1.4) follows. By Theorem 1, when deriving the small ball
asymptotics the relation (3.2) allows us to  put $\varphi$ instead
of $\cal N$ in all formulas. So, we obtain
$$L'(u)\sim -\int\limits_0^1\frac {\varphi(\lambda)\
d\lambda}{(1+2u\lambda)^2}=-\frac {\varphi(\frac 1u)}u\cdot
\int\limits_0^{\infty}\chi_{[0,u]}(t)\cdot\frac {\varphi(\frac
tu)}{\varphi(\frac 1u)}\cdot \frac {dt}{(1+2t)^2}. $$

\noindent The integrand tends pointwise to $\frac 1{(1+2t)^2}$ as
$u\to\infty$. Moreover, for any $p>0$ the function $t\to
t^p\varphi(\frac tu)$ is increasing in a neighborhood of the
origin, and hence
$$\chi_{[0,u]}(t)\cdot\frac {\varphi(\frac tu)} {\varphi(\frac
1u)} \le Ct^{-p}.$$
This provides a summable majorant, and the
Lebesgue Theorem gives
$$L'(u)\sim -\frac {\varphi(\frac 1u)}u \cdot\int
\limits_0^{\infty} \frac {dt}{(1+2t)^2}=-\frac {\varphi(\frac
1u)}{2u}, \qquad u\to\infty,
$$
and (3.4) follows.

Further,
$$L(u)\sim-\int\limits_0^u\frac {\varphi(\frac tu)\ dt}{1+2t}.
$$

Note that, given $A>0$,
$$\int\limits_0^A\frac {\varphi(\frac tu)\ dt}{1+2t}=
O(\varphi(1/u)),\qquad u\to\infty,$$

\noindent while for $u>A$
\begin{multline*}
\int\limits_A^u\frac {\varphi(\frac t u)\ dt}{1+2t}=\int \limits
_1^{u/A} \frac {\varphi(\frac {Ax}u)\ dx}{\frac 1A+2x} = \frac
12\int \limits_1^{u/A} \varphi(Ax/u)\frac {dx}x \cdot
(1+o_A(1))=\\
=\frac {\psi(\frac Au)}2\ (1+o_A(1))\sim\frac {\psi(\frac 1u)}2\
(1+o_A(1)).
\end{multline*}

Finally, since $ur\sim\frac {\varphi(\frac 1u)}2$ as $r\to0$, the
relation (3.5) gives $L(u)+ur\sim-\frac {\psi(\frac 1u)}2$, and
(3.3) follows.\hfill$\square$\medskip

As an example we consider a set of stationary Gaussian processes
${\cal R}_{C,\alpha}$ ($C,\alpha>0$), with zero mean-value and the
spectral density
$${\cal K}_{{\cal R}_{C,\alpha}}(\xi)=\exp(-C|\xi|^{\alpha}),
\qquad \xi\in\mathbb R.$$

\noindent The corresponding covariances $G_{{\cal R}_{C,\alpha}}$
are smooth functions. For example, it is well-known that
$$G_{{\cal R}_{C,1}}(s,t)=\frac C{\pi(C^2+(s-t)^2)}, \qquad
G_{{\cal R}_{C,2}}(s,t)=\frac 1{2\sqrt{\pi C}}\exp\left(-\frac
{(s-t)^2}{4C}\right).
$$

The eigenvalue asymptotics of the integral operators at the finite
interval with kernels of this type was treated in remarkable paper
\cite{W}. We underline that for $\alpha<1$ Theorem 1 \cite{W}
provides the asymptotics of $\lambda_n$, while for $\alpha\ge1$
Theorems 2 and 3 \cite{W} give only the asymptotics of
$\ln(\lambda_n)$ that enables only to obtain the asymptotics of the counting 
function. In order to provide a unified approach, we find the asymptotics of 
${\cal N}(\lambda)$ for all $\alpha$. Namely, as $\lambda\to 0$,
$${\cal N}(\lambda)\sim\varphi(\lambda)\equiv\left\{\begin{array}{rll}
\frac 1{\pi C^{\frac 1\alpha}} & \cdot\ln^{\frac 1\alpha}(\frac
1\lambda), &
\ \alpha<1;\\ \\
\frac 1{\pi {\mathfrak C}} & \cdot\ln(\frac 1\lambda), &
\ \alpha=1;\\ \\
\frac 1{2-2/\alpha} & \cdot\frac {\ln(\frac
1\lambda)}{\ln\ln(\frac 1\lambda)}, &
\ \alpha>1;\\
\end{array}
\right.$$

\noindent Here
$${\mathfrak C}=\frac {{\bf K}(\mbox{sech}(\pi/2C))}{{\bf K}
(\tanh(\pi/2C))}$$
while $\bf K$ is the complete elliptic
integral of the first kind.\medskip

Let $\alpha<1$. Then the relation (3.4) reads
$$r\sim\frac {\ln^{\frac 1\alpha}(u)}{2\pi C^{\frac
1\alpha}u} \quad \Longleftrightarrow \quad u\sim\frac {\ln^{\frac
1\alpha}(\frac 1r)} {2\pi C^{\frac 1\alpha}r}.
$$

Therefore, (3.3) provides, after substitution $r=\ep^2$,
$$\ln{\bf P}\left\{\|{\cal R}_{C,\alpha}\|\le\ep\right\} \sim-
\frac {\alpha\ln^{\frac {\alpha+1}\alpha} (u)} {2\pi(\alpha+1)
C^{\frac 1\alpha}}\sim- \left(\frac 2C\right) ^{\frac
1\alpha}\cdot \frac {\alpha\ln^{\frac {\alpha+1}\alpha} (\frac
1\ep)}{(\alpha+1)\pi},\qquad \ep\to0.\eqno(3.6)
$$

Similarly, for $\alpha=1$
$$\ln{\bf P}\left\{\|{\cal R}_{C,\alpha}\|\le\ep\right\} \sim-
\frac {\ln^2(\frac 1\ep)}{\pi{\mathfrak C}},\qquad \ep\to0.
\eqno(3.7)
$$

For $\alpha>1$ (3.4) gives
$$r\sim \frac 1{4-\frac 4\alpha}\cdot\frac {\ln(u)}{u\ln\ln(u)}
\quad \Longleftrightarrow \quad u\sim \frac 1{4-\frac 4\alpha}
\cdot\frac {\ln(\frac 1r)}{r\ln\ln(\frac 1r)},
$$

\noindent and we obtain after simple calculation
$$\ln{\bf P}\left\{\|{\cal R}_{C,\alpha}\|\le\ep\right\} \sim-
\frac 1{2-\frac 2\alpha} \cdot\frac {\ln^2(\frac 1\ep)}
{\ln\ln(\frac 1\ep)},\qquad \ep\to0. \eqno(3.8)
$$

\noindent Thus, the logarithmic small ball asymptotics in this case does not 
depend on $C$.\medskip 

{\bf Remark 3}. The processes ${\cal R}_{C,\alpha}$ were studied in 
\cite{atall} where the order of decreasing in logarithmic scale was
obtained for the small ball probabilities in $\sup$-norm. We note 
that our approach gives an alternative proof of the key upper estimate 
in \cite{atall} due to a trivial relation $\|X\|\le \sup|X|$.

\bigskip
I am grateful to Prof. M.A.~Lifshits for useful discussions and to Prof. 
M.Z.~Solomyak for the hint to the reference \cite{W}. The most part of this 
paper was written during my visit to Link\"oping University in Sweden. I 
would like to thank Prof. V.A.~Kozlov for excellent working conditions.


\bigskip

\begin{thebibliography}{AFTL2}

\bibitem[AILZ]{atall} 
F.~Aurzada, I.A.~Ibragimov, M.A.~Lifshits, H.~van~ Zanten, {\em Small
deviations of smooth stationary Gaussian processes}, available at
{\tt www.arXiv.org/abs/0803.4238}; submitted to Theor. Probab. Appl.

\bibitem[GHLT1]{GHLT1} 
F.~Gao, J.~Hannig, T.-Y.~Lee, F.~Torcaso, {\em Laplace
transforms via Hadamard factorization with applications to small
ball probabilities,} Electronic J. Probab. {\bf8(13)} (2003), 1--20.

\bibitem[GHLT2]{GHLT}
F.~Gao, J.~Hannig, T.-Y.~Lee, F.~Torcaso,
{\em Exact $L^2$-small balls of Gaussian processes}, J.  Theor.
Probab., {\bf 17} (2004), 503--520.

\bibitem[KNN]{KNN} 
A.I.~Karol, A.I.~Nazarov, Ya.Yu.~Nikitin, {\em Small ball
probabilities for Gaussian random fields and tensor products of
compact operators}, Trans. Amer. Math. Soc., {\bf360(3)} (2008), 1443--1474.

\bibitem[Li]{Li}
W.V.~Li, {\em Comparison results for the lower tail of
Gaussian semi\-norms}, J. Theor. Probab., {\bf 5} (1992), 1--31.

\bibitem[LS]{LS}
W.V.~Li, Q.M.~Shao, {\em Gaussian processes:
inequalities, small ball probabilities and applications},
Stochastic Processes: Theory and Methods. Handbook of Statistics,
{\bf 19} (2001), C.R.Rao and D.Shanbhag (Eds), 533--597.

\bibitem[Lf1]{Lf2} 
M.A.~Lifshits, {\em On the lower tail probabilities of some
random series}, Ann. Probab., {\bf25(1)} (1997), 424--442.

\bibitem[Lf2]{Lf}
M.A.~Lifshits, {\em Asymptotic behavior of small ball
probabilities}, Prob.\ Theory and Math.\ Stat., 1999.
B.Grigelionis et al.\ (Eds), Proc. VII International Vilnius
Conference (1998), VSP/TEV, 453--468.

\bibitem[Lf3]{site}
M.A.~Lifshits, {\em Bibliography on small deviation probabilities}.
Compilation available at
{\tt www.proba.jussieu.fr/pageperso/smalldev/biblio.html}

\bibitem[Na1]{Na2}
A.I.~Nazarov,
{\em On the sharp constant in the small ball asymptotics of some
Gaussian processes under $L_2$-norm}, Probl. Mat. Anal. {\bf26} (2003),
179--214 (Russian); English transl.: J. Math. Sci. {\bf117(3)} (2003),
4185--4210.

\bibitem[Na2]{Na3} 
A.I.~Nazarov, {\em Logarithmic asymptotics of small
deviations for some Gaussian processes in the ${L_2}$-norm with
respect to a self-similar measure}, ZNS POMI, {\bf311} (2004),
190--213 (Russian); English transl.: J. Math. Sci., {\bf133(3)}
(2006), 1314--1327.

\bibitem[Na3]{Na1}
A.I.~Nazarov, {\em Exact $L_2$-small ball asymptotics of
Gaussian processes and the spectrum of boundary value problems
with "non-separated"\ boundary conditions}, available at
{\tt http://arxiv.org/abs/0710.1408}; to appear in J. Theor. Probab.

\bibitem[NN1]{NN1}
A.I.~Nazarov, Ya.Yu.~Nikitin, {\em Exact $L_2$-small
ball behavior of integrated Gaussian processes and spectral
asymptotics of boundary value problems}, Probab. Theory 
Rel. Fields, {\bf 129} (2004), 469--494.

\bibitem[NN2]{NN2} 
A.I.~Nazarov, Ya.Yu.~Nikitin, {\em Logarithmic ${L_2}$-small
ball asymptotics for some fractional Gaussian processes}, Teor.
Ver. Primen., {\bf 49} (2004), N4, 695--711 (Russian); English
transl.: Theor. Probab. Appl., {\bf 49} (2005), N4, 645--658.

\bibitem[NP]{NP} 
A.I.~Nazarov, R.S.~Pusev, {\em Exact $L_2$-small ball
asymptotics for some weighted Gaussian processes,} Preprint of St.Petersburg 
Math. Soc. N 2006-1 (Russian); available at 
{\tt http://www.mathsoc.spb.ru/preprint/2006/index.html\#01};
to appear in ZNS POMI.

\bibitem[NSh]{NSh} 
A.I.~Nazarov, I.A.~Sheipak, {\em Small deviations of
Gaussian processes in $L_2$-norm with respect to a degenerate
self-similar measure}, in preparation.

\bibitem[Se]{Se} 
E.~Seneta, {\em Regularly Varying Functions}, Lect. Notes Math., 
{\bf 508} (1976).

\bibitem[S]{S}
G.N.~Sytaya, {\em On some asymptotic representations of the
Gaussian measure in a Hilbert space}, Theory of Stochastic
Processes, Kiev, {\bf 2} (1974), 93--104 (Russian).

\bibitem[W]{W} 
H.~Widom, {\em Asymptotic behavior of the eigenvalues of
certain integral equations II}, Arch. Rat. Mech. Anal., {\bf17(3)}
(1964), 215--229.


\end{thebibliography}
\end{document}